\documentclass[11pt]{article}

\usepackage{graphics}
\usepackage{graphicx}
\usepackage{color}
\usepackage[all]{xy}
\newfont{\bbb}{msbm10 scaled\magstephalf}
\newfont{\sbbb}{msbm7 scaled\magstephalf}

\def\C{\mbox{\bbb{C}}}

\def\R{\mbox{\bbb{R}}}
\def\Z{\mbox{\bbb{Z}}}
\def\cd{\C^d}
\def\rd{\R^d}
\def\rk{\R^k}
\def\rn{\R^n}
\def\rtwo{\R^2}
\def\zd{\Z^d}
\def\td{T^d}
\def\rddu{(\rd)^*}
\def\rtwodu{(\R^2)^*}
\def\rndu{(\rn)^*}

\def\et1{e^{2\pi i\theta_1}}
\def\etd{e^{2\pi i\theta_d}}

\def\vz{\underline{z}}
\def\zjs{|z_j|^2}

\def\cl{{\mathcal C}_\lambda}
\def\c0{{\mathcal C}_0}
\def\D{\Delta}
\def\Dsh{\D^{\sharp}}
\def\xd{X_1,\ldots,X_d}
\def\ld{\lambda_1,\ldots,\lambda_d}

\def\G{\Gamma}
\def\Gsh{\G^{\sharp}}

\def\gsh{\gamma^{\sharp}}
\def\vt{\tilde{V}}
\def\ut{\tilde{U}}
\def\fia{\tau_{\alpha}}
\def\ga{\Gamma_{\alpha}}
\def\va{V_{\alpha}}
\def\vta{\vt_{\alpha}}
\def\fib{\tau_{\beta}}
\def\gb{\G_{\beta}}
\def\vb{V_{\beta}}
\def\vtb{\vt_{\beta}}
\def\gab{g_{\alpha\beta}}

\def\wt{\tilde{W}}
\def\vsh{V^{\sharp}}
\def\wsh{W^{\sharp}}
\def\fbar{\bar{f}}
\def\fsh{f^{\sharp}}

\newfont{\frak}{eufm10 scaled\magstep1}
\newfont{\sfrak}{eufm8 scaled\magstep1}

\def\n{\mbox{\frak n}}
\def\l{\mbox{\frak l}}

\newcommand{\proof}{\mbox{\bf Proof.\ \ }}

\def\squareforqed{\hbox{\rlap{$\sqcap$}$\sqcup$}}
\def\qed{\ifmmode\else\unskip\quad\fi\squareforqed}
\def\smartqed{\def\qed{\ifmmode\squareforqed\else{\unskip\nobreak\hfil
\penalty50\hskip1em\null\nobreak\hfil\squareforqed
\parfillskip=0pt\finalhyphendemerits=0\endgraf}\fi}}

\newtheorem{thm}{Theorem}[section]
\newtheorem{defn}[thm]{Definition}
\newtheorem{prop}[thm]{Proposition}
\newtheorem{remark}[thm]{Remark}
\newtheorem{lemma}[thm]{Lemma}
\title{The Symplectic Geometry of Penrose Rhombus Tilings}
\author{\sc Fiammetta Battaglia and Elisa Prato}
\date{}
\begin{document}
\maketitle
\pagestyle{myheadings}
\markboth{\sc Fiammetta Battaglia and Elisa Prato}{\sc The Symplectic
Geometry of Penrose Rhombus Tilings}
\begin{abstract}
The purpose of this article is to view Penrose rhombus tilings from
the perspective of symplectic geometry. We show that each thick
rhombus in such a tiling can be naturally associated to a highly
singular $4$--dimensional compact symplectic space $M_R$, while each
thin rhombus can be associated to another such space $M_r$; both
spaces are invariant under the Hamiltonian action of a
$2$--dimensional quasitorus, and the images of the corresponding
moment mappings give the rhombuses back. The spaces $M_R$ and $M_r$
are diffeomorphic but not symplectomorphic.
\end{abstract}
\centerline{
{\small {\em Mathematics Subject Classification}. Primary: 53D20. 
Secondary: 52C23}}
\section*{Introduction}
We start by considering a Penrose tiling by thick and thin rhombuses
(cf. \cite{pen1}). Rhombuses are very special examples of simple
convex polytopes. Because of the Atiyah, Guillemin--Sternberg
convexity theorem \cite{a,gs}, convex polytopes can arise as images
of the moment mapping for Hamiltonian torus actions on compact
symplectic manifolds. For example, simple convex polytopes that are
rational with respect to a lattice $L$ and satisfy an additional
integrality condition, correspond to symplectic toric manifolds.
More precisely, the Delzant theorem \cite{d} tells us that to each
such polytope in $\rndu$, there corresponds a compact symplectic
$2n$--dimensional manifold $M$, endowed with the effective
Hamiltonian action of a torus of dimension $n$. As it turns out, the
polytope is exactly the image of the corresponding moment mapping.
One of the striking features of Delzant's theorem is that it gives
an explicit procedure to obtain the manifold corresponding to each
given polytope as a symplectic reduced space. This correspondence
may be applied to each of the rhombuses in a Penrose tiling
separately. However, the rhombuses in a Penrose tiling, though
simple and convex, are not simultaneously rational with respect to
the same lattice. Therefore we cannot apply the Delzant procedure
simultaneously to all rhombuses in the tiling. However, if we
replace the lattice with a quasilattice (the $\Z$--span of a set of
$\R$--spanning vectors) and the manifold with a suitably singular
space, then it is possible to apply a generalization of the Delzant
procedure to arbitrary simple convex polytopes that was given by the
second named author in \cite{p1}. According to this result, to each
simple convex polytope in $\rndu$, and to each suitably chosen
quasilattice $Q$, one can associate a family of compact symplectic
$2n$--dimensional {\em quasifolds} $M$, each endowed with the
effective Hamiltonian action of the quasitorus $\rndu / Q$, having
the property that the image of the corresponding moment mapping is
the polytope itself. Quasifolds are generalizations of manifolds and
orbifolds that were introduced in \cite{p1}. A local model for a
$k$--dimensional quasifold is given by the topological quotient of
an open subset of $\rk$ by the action of a finitely generated group.
A $k$--dimensional quasifold is a topological space admitting an
atlas of $k$--dimensional local models that are suitably glued
together. It is a highly singular, usually not even Hausdorff,
space. A quasitorus is the natural generalization of a torus in this
setting.

We remark that, unlike the Delzant case, we do not have here a
one--to--one correspondence between polytopes and symplectic spaces.
In fact, there is much more freedom of choice in the generalized
Delzant construction, and infinitely many symplectic quasifolds will
map to the same polytope. More precisely, given any suitable
quasilattice $Q$, the construction will yield one quasifold for each
choice of a set of vectors $X_1,\ldots,X_d$ in $Q$, each of which is
orthogonal and pointing inwards to one of the $d$ different facets
of the polytope. The striking fact in the case of a Penrose tiling
is that there is a natural choice of a quasilattice and of a set of
inward--pointing vectors, and therefore a natural choice of a
privileged quasifold mapping to each tile. In fact, let us consider
a Penrose rhombus tiling where all rhombuses have edges of length
$1$; we will see in Section~1 that such a tiling determines a star
of five unit vectors, pointing to the vertices of a regular
pentagon. One of the features of the tiling is that the four edges
of any given rhombus in the tiling are orthogonal to two of these
five vectors. These two vectors, and their opposites, will be our
natural choice of the four inward--pointing vectors in the Delzant
construction. The span over the integers of the five vectors in the
star is dense in $\R^2$, this is the quasilattice that is associated
to the tiling, and will be our choice of quasilattice.

Moreover, we show that all the different rhombuses yield only two
possible compact symplectic quasifolds: one for all the thick
rhombuses, $M_R$, and one for all the thin ones, $M_r$. Both are
global quotients of a product of two spheres modulo the action of a
finitely generated group; they are diffeomorphic, but not
symplectomorphic. Notice that in general a quasifold is not a global
quotient of a manifold modulo the action of a finitely generated
group, in fact most of the examples are not, see the examples worked
out in \cite{p1} and the quasifold corresponding to the Penrose kite
in \cite{kite}.

We remark that quasilattices are quasiperiodic structures underlying
quasiperiodic tilings and the atomic order of quasicrystals, which
were discovered in the eighties  by observing that the diffraction
pattern of such materials is not periodic, but quasiperiodic, with
5--fold rotational symmetries (cf. \cite{quasicristalli}). In
forthcoming work \cite{3d} we give a symplectic interpretation of
three--dimensional analogues of Penrose tilings, that represent
structure models for icosahedral quasicrystals.

The paper is structured as follows. In the first section we recall
some basic facts on Penrose rhombus tilings. In the second section
we recall from \cite {p1} the notions of symplectic quasifold, of
quasitorus actions, and the generalization of the Delzant procedure.
In the third section we apply all of the above to construct the
quasifolds $M_R$ and $M_r$, from special choices of one thick
rhombus and of one thin rhombus in a given Penrose tiling. Finally
in the fourth section we show that all the other thick rhombuses of
the tiling correspond to $M_R$, that all the other thin rhombuses
correspond to $M_r$ and that $M_R$ and $M_r$ are diffeomorphic but
not symplectomorphic.
\section{Penrose Rhombus Tilings}
\subsection{Penrose Rhombuses}
\begin{figure}
\begin{center}
\includegraphics{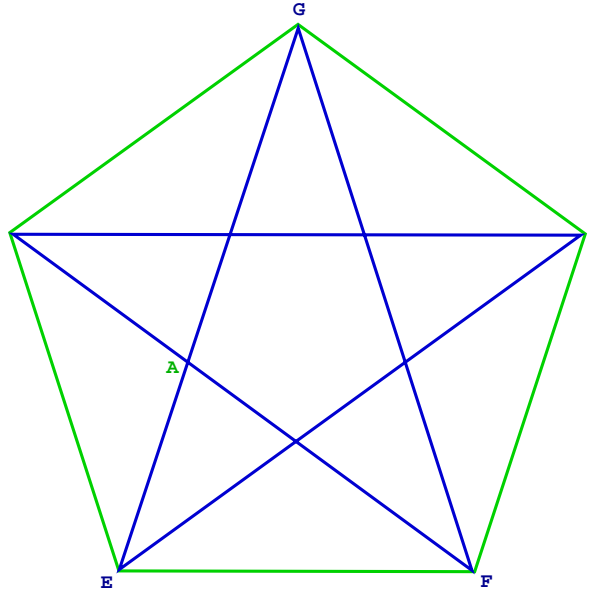}
\end{center}
\caption{The pentagram} \label{pentagram}
\end{figure}
\begin{figure}
\begin{center}
\includegraphics{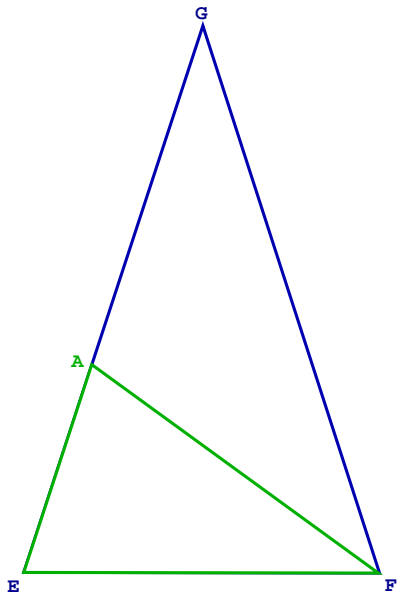}
\end{center}
\caption{The golden triangle and its decomposition}
\label{goldentriangle}
\end{figure}
\begin{figure}
\begin{center}
\includegraphics{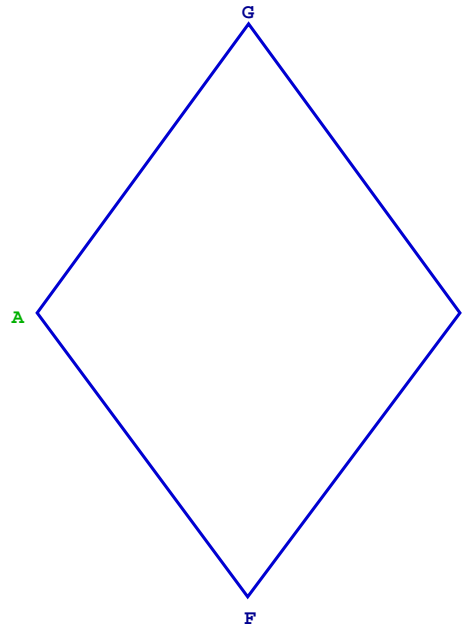}
\end{center}
\caption{The thick rhombus} \label{thickrhombus}
\end{figure}
\begin{figure}
\begin{center}
\includegraphics{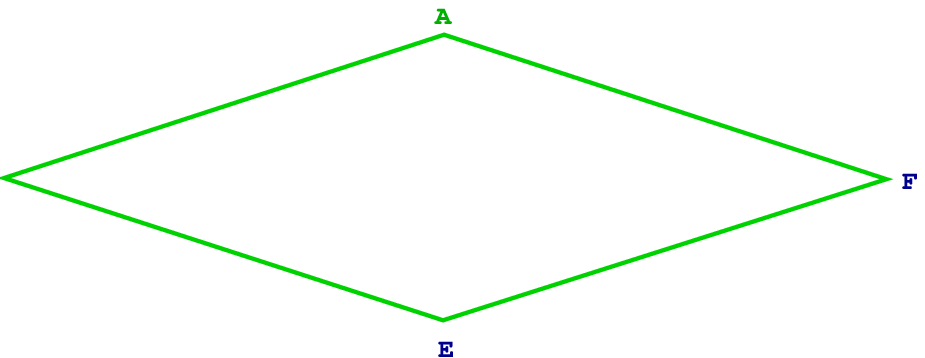}
\end{center}
\caption{The thin rhombus} \label{thinrhombus}
\end{figure}
Let us now recall the procedure for obtaining the Penrose rhombuses
from the pentagram. For a proof of the facts that are needed we
refer the reader to \cite{l}, and for additional historical remarks
we refer the reader to \cite{p2}. Let us consider a regular pentagon
whose edges have length one and let us consider the corresponding
inscribed pentagram, as in Figure~\ref{pentagram}. It can be shown
that the ratio of the diagonal to the side of the pentagon is equal
to the {\em golden ratio},
$\phi=\frac{1+\sqrt{5}}{2}=2\cos{\frac{\pi}{5}}$. Therefore the
triangle having vertices $\mbox{E, F, G}$ is a {\em golden
triangle}, which is, by definition, an isosceles triangle with a
ratio of side to base given by $\phi$. This triangle decomposes into
the two smaller triangles of vertices $\mbox{E, F, A}$ and $\mbox{F,
G, A}$, respectively (see Figure~\ref{goldentriangle}). The first
one is itself a golden triangle. Using the fundamental relation
\begin{equation}\label{phi}
\phi=1+\frac{1}{\phi}
\end{equation}
one can show that the second one is a {\em golden gnomon}, which is,
by definition, an isosceles triangle with a ratio of side to base
given by $\frac{1}{\phi}$. Now, if we consider the union of the
golden gnomon with its reflection with respect to the
$\mbox{FG}$--axis, we obtain the {\em thick rhombus} (see
Figure~\ref{thickrhombus}), while to obtain the {\em thin rhombus}
we consider the union of the smaller golden triangle with its
reflection with respect to the $\mbox{EG}$--axis (see
Figure~\ref{thinrhombus}). Notice that the angles of the thick
rhombus  measure $2\pi/5$ and $3\pi/5$, while the angles of the thin
rhombus measure $\pi/5$ and $4\pi/5$.
\subsection{The Tiling Construction}\label{penrose}
\begin{figure}
\begin{center}
\includegraphics{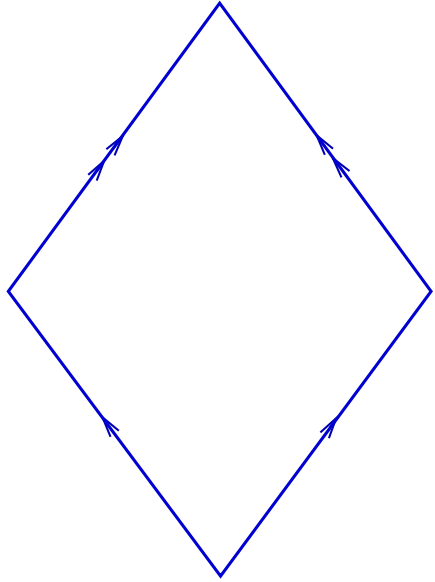}
\end{center}
\caption{Matching rules for the thick rhombus}
\label{matchingrules1}
\end{figure}

\begin{figure}
\begin{center}
\includegraphics{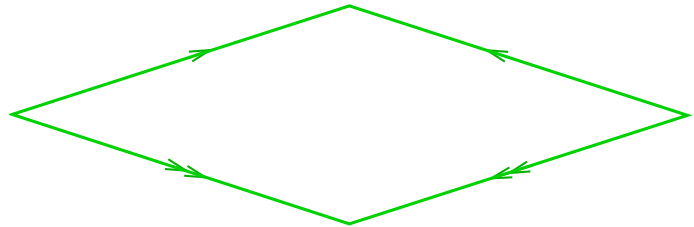}
\end{center}
\caption{Matching rules for the thin rhombus} \label{matchingrules2}
\end{figure}

\begin{figure}
\begin{center}
\includegraphics[width=0.5\linewidth]{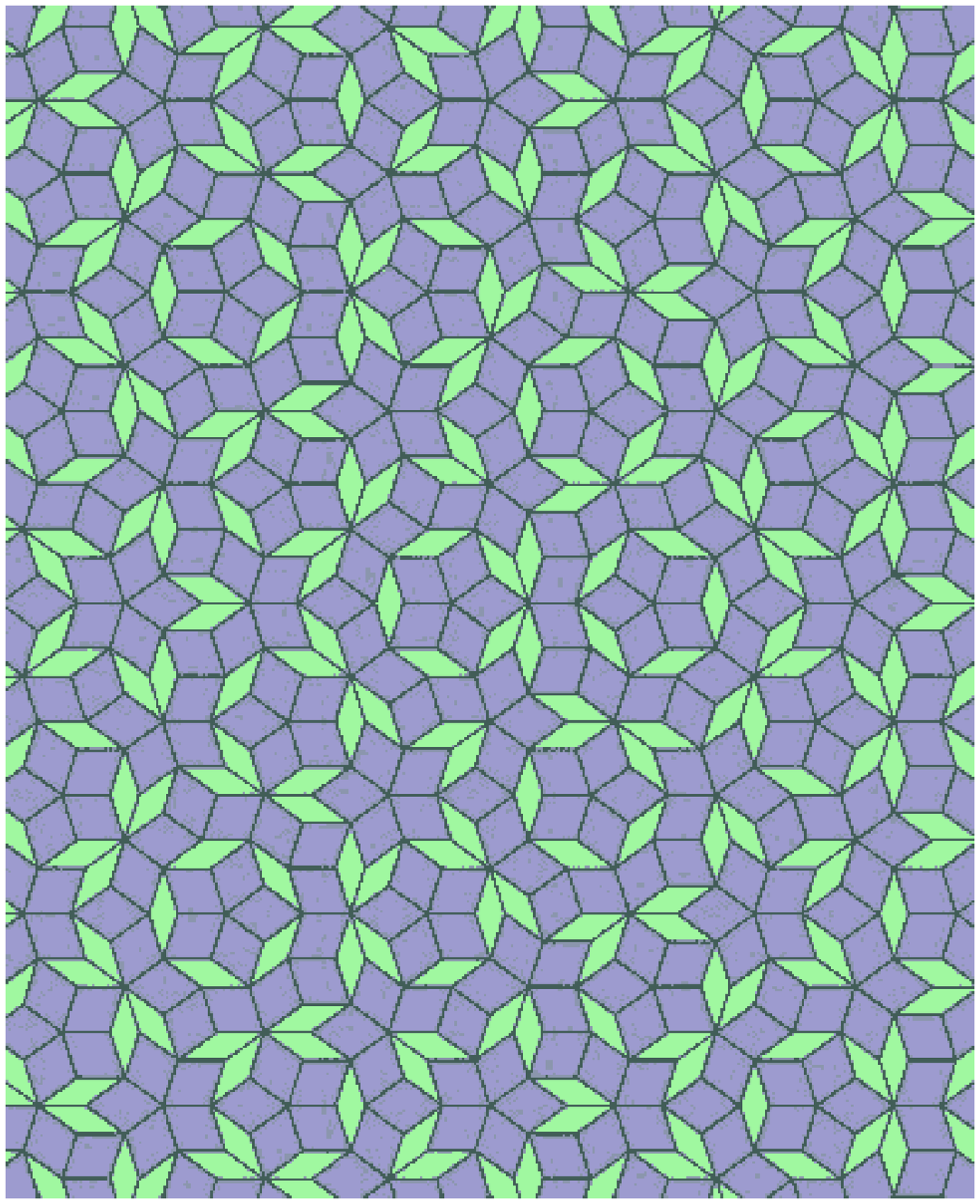}
\end{center}
\caption{A Penrose rhombus tiling. Figure by D. Austin [2],
reprinted courtesy of the AMS} \label{penrosetiling}
\end{figure}

We now briefly recall some basic facts about Penrose rhombus
tilings; for a deeper analysis of this important subject we refer
the reader to the original paper by Penrose \cite{pen1}, to his
subsequent works \cite{pen2,pen3}, to Austin's articles
\cite{ams1,ams2} and finally, for a review, to the books by
Gr\"unbaum and Shephard \cite{grun} and by Senechal \cite{S}. We
start with the following
\begin{defn}[Penrose rhombus tiling] A Penrose rhombus tiling is a
tiling of the plane by thick and thin rhombuses that obey the
matching rules shown in Figures~\ref{matchingrules1} and
\ref{matchingrules2}.
\end{defn}
Consider a Penrose tiling in $(\R^2)^*$ with rhombuses having edges
of length $1$ (Figure~\ref{penrosetiling});
it is well--known that there are uncountably many such tilings and
that each of them is non periodic. The key remark in our set--up is
that there exist a quasilattice, and a set of vectors contained in
it, that are naturally associated to any such tiling. First of all
we need to give the formal definition of quasilattice:
\begin{defn}[Quasilattice]
Let $V$ be a real vector space. A {\em quasilattice} in $V$ is the
Span over $\Z$ of a set of $\R$--spanning vectors $V_1,\ldots,V_d$
of $V$.
\end{defn}
Notice that $\hbox{Span}_{\Z}\{V_1,\dots,V_d\}$ is a lattice if and
only if it admits a set of generators which is a basis of $V$.
Consider now, in $(\R^2)^*$, the star ${\mathcal S}^*$ of five unit
vectors
$$\left\{\begin{array}{l}
Y^*_0=(0,1)\\
Y^*_1=\frac{1}{2}(-\sqrt{2+\phi},\frac{1}{\phi})\\
Y^*_2=\frac{1}{2}(-\frac{1}{\phi}\sqrt{2+\phi},-\phi,)\\
Y^*_3=\frac{1}{2}(\frac{1}{\phi}\sqrt{2+\phi},-\phi)\\
Y^*_4=\frac{1}{2}(\sqrt{2+\phi},\frac{1}{\phi})
\end{array}
\right.$$
\begin{figure}
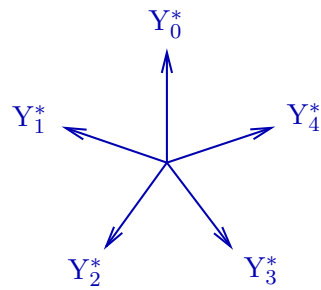

\begin{center}
\input stargridduale.pstex_t
\end{center}
\caption{Dual star ${\mathcal S}^*$} \label{stargriddual}
\end{figure}
Let $R$ be the quasilattice generated by the vectors of ${\mathcal
S}^*$, namely $$R=\hbox{Span}_{\Z}\{Y^*_0,\dots,Y^*_4\}.$$ The
quasilattice $R$ is not a lattice, it is dense in $\R^2$ and a
minimal set of generators of $R$ is made of $4$ vectors. The
following statement describes the relationship between the Penrose
rhombus tilings considered and the quasilattice $R$, together with
its generators.

Let us consider the five thick rhombuses that are determined by
$Y_k^*,Y_{k+1}^*$, for $k=0,\dots,4$, we denote them by $\D_R^k$,
and the five thin rhombuses that are determined by
$Y_k^*,Y^*_{k+2}$, for $k=0,\dots,4$, we denote them by $\D_r^k$ (we
are assuming here $Y^*_5=Y^*_0$ and $Y^*_6=Y^*_1$).

Let us consider any Penrose rhombus tiling $\mathcal T$ with
rhombuses having edges of length $1$ and denote one of its edges by
$AB$. From now on we will choose our coordinates so that $A=O$ and
$B-A=Y^*_0$.
\begin{prop}\label{rotazioni} Let $\mathcal T$ be a Penrose rhombus tiling with
rhombuses having edges of length $1$. Then each rhombus is the
translate of either a thick rhombus $\D_R^k$, $k=0,\ldots,4$, or of
a thin rhombus $\D_r^k$, $k=0,\ldots,4$. Moreover each vertex of the
tiling lies in the quasilattice $R$.\end{prop} \proof The argument
is very simple. Let $C$ be a vertex of the tiling that is different
from $0$ and the above vertex $B$. We can walk from $B$ to $C$ on a
path made of subsequent edges of the tiling. We denote the vertices
of the broken line thus obtained by
$T_0=A,T_1=B,\dots,T_j,\dots,T_m=C$. The angle of the broken line at
each vertex $T_j$ is necessarily a multiple of $\pi/5$. Therefore
each vector $V_j=T_j-T_{j-1}$ is one of the vectors $\pm Y^*_k$,
$k=0,\ldots,4$. Since $C-A=T_m-T_0=V_m+\cdots+V_1$ the vertex $C$
lies in $R $ and each rhombus having $C$ as vertex has edges
parallel to two vectors of $S^*$.\qed

Consider now the star of vectors $\mathcal S$ in $\R^2$ given by:
\begin{equation}\label{star}\left\{\begin{array}{l}
Y_0=(1,0)\\
Y_1=(\cos{\frac{2\pi}{5}},\sin{\frac{2\pi}{5}})=\frac{1}{2}(\frac{1}{\phi},\sqrt{2+\phi})\\
Y_2=(\cos{\frac{4\pi}{5}},\sin{\frac{4\pi}{5}})=\frac{1}{2}(-\phi,\frac{1}{\phi}\sqrt{2+\phi})\\
Y_3=(\cos{\frac{6\pi}{5}},\sin{\frac{6\pi}{5}})=\frac{1}{2}(-\phi,-\frac{1}{\phi}\sqrt{2+\phi})\\
Y_4=(\cos{\frac{8\pi}{5}},\sin{\frac{8\pi}{5}})=\frac{1}{2}(\frac{1}{\phi},-\sqrt{2+\phi})
\end{array}
\right.\end{equation}
\begin{figure}
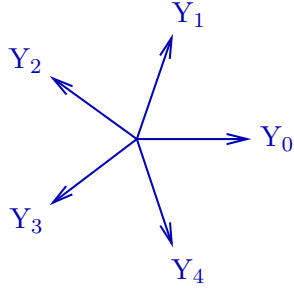

\begin{center}
\input stargrid.pstex_t
\end{center}
\caption{Star $\mathcal S$} \label{stargrid}
\end{figure}
It is easy to check that for each $k$ the four vectors $\pm Y_k,\pm
Y_{k+1}$ are each orthogonal and inward--pointing to one of the
different edges of the thick rhombus $\D_R^k$. In the same way the
four vectors $\pm Y_k,\pm Y_{k+2}$ are each orthogonal and
inward--pointing to one of the different edges of the thin rhombus
$\D_r^k$. By Proposition~\ref{rotazioni}, the same is true for each
thick and thin rhombus of any given tiling.

We denote by $Q$ the quasilattice generated by the vectors of $S$,
namely
\begin{equation}\label{qu}
Q=\hbox{Span}_{\Z}\{Y_0,\dots,Y_4\}.\end{equation}

The following relations are necessary for determining the groups
involved in the construction of the quasifolds corresponding to the
rhombuses. If we write each $Y_k$, in complex notation, as
$Y_k=e^{\frac{2 k \pi i}{5}}$, then it can be easily verified that
\begin{equation}\label{rombone}\begin{array}{l}
Y_{k+2}=-Y_{k}+\frac{1}{\phi}Y_{k+1}\\
Y_{k+4}=\frac{1}{\phi}Y_{k}-Y_{k+1}\\
\end{array}\end{equation}
and
\begin{equation}\label{rombino}\begin{array}{l}
Y_{k+3}=-Y_{k+2}-\phi Y_k\\
Y_{k+4}=-\phi Y_{k+2}-Y_k.\\
\end{array}\end{equation}
Moreover $Y_0=-(Y_1+Y_2+Y_3+Y_4)$
and $\{Y_1,Y_2,Y_3,Y_4\}$ is a minimal set of generators of $Q$.
\section{Symplectic Quasifolds}
Let us recall the definition of quasifold; we refer to the article
\cite{p1} for the missing details and proofs. We begin by defining
the
\begin{defn}[Quasifold model]\label{model}{\rm
    Let $\vt$ be a connected open subset of
    $\rk$ and let $\G$ be a finitely generated group acting
    smoothly on $\vt$
    so that the set of points, $\vt_0$, where the action is free,
    is connected and dense. Consider the space of orbits, $\vt/\G$,
    of the action of the group
    $\G$ on $\vt$, endowed with the quotient topology, and the
    canonical projection $p\;\colon\;\vt\rightarrow \vt/\G$.
    A {\em quasifold model} of dimension $k$ is the triple $(\vt/\G,p,\vt)$,
    shortly denoted
    $\vt/\G$.}
\end{defn}
\begin{defn}[Submodel]{\rm
Consider a model $(\vt/\G,p,\vt)$ and let $W$ be an open subset of
$\vt/\G$. We will say that $W$ is a submodel of $(\vt/\G,p,\vt)$, if
$(W,p,p^{-1}(W))$ defines a model.}
\end{defn}
\begin{remark}\label{simplyc}{\rm Consider a model of dimension
$k$, $(\vt/\G,p,\vt)$, such that there exists a covering
$\pi\,\colon\,\vsh\rightarrow\vt$, where $\vsh$ is an open subset of
$\rk$ acted on by a finitely generated group $\Pi$ in a smooth, free
and proper fashion with $\vt=\vsh/\Pi$. Consider the extension of
the group $\G$ by the group $\Pi$
$$1\longrightarrow \Pi\longrightarrow\Gsh
\longrightarrow\G\longrightarrow 1$$ defined as follows
$$\Gsh=\left\{\;\gsh\in\mbox{Diff}(\vsh)\;|\;\exists\;
\gamma\in\Gamma\;\mbox{s. t.}\;\pi(\gsh(u^{\sharp}))=\gamma\circ
\pi(u^{\sharp})\;\forall\; u^{\sharp}\in\vsh\;\right\}.$$ It is easy
to verify that $\Gsh$ is finitely generated, that it acts on $\vsh$
according to the assumptions of Definition~\ref{model} and that
$\vt/\G=\vsh/\Gsh$. Let $p^{\sharp}=p\circ\pi$, we will then say
that the model $(\vsh/\Gsh,p^{\sharp},\vsh)$ is a covering of the
model $(\vt/\G,p,\vt)$. }\end{remark}
\begin{defn}[Smooth mapping, diffeomorphism of models]\label{diffeomodelli}{\rm
Given two models $(\vt/\G, p, \vt)$ and $(\wt/\D,q,\wt)$, a mapping
$f\,\colon\,\vt/\G\longrightarrow\wt/\D$ is said to be {\em smooth}
if there exist coverings in the sense of Remark~\ref{simplyc} of the
two given models, $(\vsh/\Gsh,p^{\sharp},\vsh)$ and $(\wsh/\Dsh,
q^{\sharp},\wsh)$, and  a smooth mapping
$\fsh\,\colon\,\vsh\longrightarrow\wsh$ such that $ q^{\sharp}\circ
\fsh= f\circ p^{\sharp}$; we will then say that $\fsh$ is a {\em
lift} of $f$. We will say that the smooth mapping $f$ is a {\em
diffeomorphism of models} if it is bijective and if the lift $\fsh$
is a diffeomorphism.}
\end{defn}
If the mapping $\fsh$ is a lift of a smooth mapping of models
$f\,\colon\,\ut/\G\longrightarrow \vt/\D$ so are the mappings
${\fsh}^{\gamma}(-)=\fsh(\gamma\cdot -)$, for all elements $\gamma$ in
$\Gsh$ and $^{\delta}\fsh(-)=\delta\cdot\fsh(-)$, for all elements
$\delta$ in $\Dsh$. We recall that if the mapping $f$ is a
diffeomorphism, then these are the only other possible lifts:
\begin{lemma}\label{orange}
Consider two models, $\ut/\G$ and $\vt/\D$, and let
$f\,\colon\,\ut/\G\longrightarrow\vt/\D$ be a diffeomorphism of
models. For any two lifts, $\fsh$ and $\fbar$, of the diffeomorphism
$f$, there exists a unique element $\delta$ in $\Dsh$ such that
$\fbar={}^\delta\fsh$.
\end{lemma}
\begin{lemma}\label{green} Consider two models, $\ut/\G$ and $\vt/\D$, and
a diffeomorphism $f\,\colon\,\ut/\G\longrightarrow\vt/\D$. Then, for
a given lift, $\fsh$, of the diffeomorphism $f$, there exists a group
isomorphism $F\,\colon\,\Gsh\longrightarrow\Dsh$ such that
${\fsh}^{\gamma}={}^{F(\gamma)}\fsh$, for all elements $\gamma$ in $\Gsh$.
\end{lemma}
Similarly to the notion of smooth mapping it is possible to define
other geometric objects on models, such as differential forms,
symplectic forms and vector fields.
\begin{defn}[Quasifold]
  {\rm A dimension $k$ {\em quasifold structure} on a topological space
    $M$ is the assignment of an {\em atlas}, or collection of {\em charts},
   ${\mathcal A}= \{\,(\va,\fia,\vta/\ga)\,|\,\alpha\in A\,\}$ having the following properties:
\begin{enumerate}
\item The collection $\{\,\va\,|\,\alpha\in A\,\}$ is a cover of $M$.
\item For each index $\alpha$ in $A$, the set $\va$ is open, the space
$\vta/\ga$ defines a model, and the
mapping $\fia$ is a homeomorphism of the space $\vta/\ga$ onto the
set $\va$.
\item For all indices $\alpha, \beta$ in $A$ such that
  $\va\cap\vb\neq\emptyset$, the sets $\fia^{-1}(\va\cap\vb)$ and
  $\fib^{-1}(\va\cap\vb)$ are submodels of $\vta/\ga$ and $\vtb/\gb$
  respectively and the mapping
  $$\gab=\fib^{-1}\circ\fia\,\colon\fia^{-1}(\va\cap\vb)
  \longrightarrow\fib^{-1}(\va\cap\vb)$$
  is a diffeomorphism of models. We will then say that the mapping
  $\gab$ is a {\em change of charts} and that the corresponding charts are
  {\em compatible}.
\item The atlas $\mathcal A$ is maximal, that is: if the triple
$(V,\tau,\vt/\G)$ satisfies property 2. and is compatible with all
the charts in $\mathcal A$, then $(V,\tau,\vt/\G)$ belongs to
$\mathcal A$.
\end{enumerate}
We will say that a space $M$ with a quasifold structure is a {\em
quasifold}.}
\end{defn}
\begin{remark}\label{gruppi}{\rm
Remark that, by the definition of diffeomorphism, finitely generated
groups corresponding to different charts need not be isomorphic (see
the fundamental example of the quasisphere in \cite{p1}).}
\end{remark}
\begin{remark}\label{gruppino}{\rm
To each point $m\in M$ there corresponds a finitely generated group $\G_m$
defined as follows: take a chart $(\va,\fia,\vta/\ga)$ around $m$,
then $\G_m$ is the isotropy group of $\ga$ at any point
$\tilde{v}\in \vt$ which projects down to $m$. One can check that
this definition does not depend on the choice of the chart. If all
the $\G_m$'s are finite $M$ is an orbifold, if they are trivial then
$M$ is a manifold.}\end{remark} It is possible to define on any
quasifold $M$ the notions of smooth mapping, diffeomorphism,
differential form, symplectic form and smooth vector field.
\begin{defn}[Quasitorus]
A {\em quasitorus} of dimension $n$ is the quotient $\rn/Q$, where
$Q$ is a quasilattice in $\rn$.
\end{defn}
We remark that a quasitorus is an example of quasifold covered by
one chart. At this point one can define the notion of Hamiltonian
action of a quasitorus on a symplectic quasifold, and the
corresponding moment mapping.
\section{The Tiles from a Symplectic Viewpoint}
We now outline the generalization of the Delzant procedure \cite{d}
to nonrational simple convex polytopes that is proven in \cite{p1}.
We begin by recalling what is a
\begin{defn}[Simple polytope] A dimension $n$ convex polytope $\D\subset\rndu$ is said to be {\em
simple} if there are exactly $n$ edges stemming from each
vertex.\end{defn}
Let us now consider a dimension $n$ convex
polytope $\D\subset\rndu$. If $d$ is the number of facets of $\D$,
then there exist elements $\xd$ in $\rn$ and $\ld$ in $\R$ such that
\begin{equation}\label{polydecomp}
\D=\bigcap_{j=1}^d\{\;\mu\in\rndu\;|\;\langle\mu,X_j\rangle\geq\lambda_j\;\}.
\end{equation}
\begin{defn}[$Q$--rational polytope] Let $Q$ be a quasilattice in
$\rn$. A convex polytope $\D\subset\rndu$ is said to be
$Q$--rational, if the vectors $\xd$ can be chosen in $Q$.
\end{defn}
All polytopes in $\rndu$ are $Q$--rational with respect to some
quasilattice $Q$; it is enough to consider the quasilattice that is
generated by the elements $\xd$ in (\ref{polydecomp}). Notice that
if the quasilattice is a honest lattice then the polytope is
rational.

In our situation we only need to consider the special case of simple
convex polytopes in $2$--dimensional space. Let $Q$ be a
quasilattice in $\rtwo$ and let $\D$ be a simple convex polytope in
the space $\rtwodu$ that is $Q$--rational. Consider the space $\cd$
endowed with the standard symplectic form $\omega_0=\frac{1}{2\pi
i}\sum_{j=1}^d dz_j\wedge d\bar{z}_j$ and the standard action of the
torus $\td=\rd/\zd$:
$$
\begin{array}{cccccl}
\tau\,\colon& \td&\times&\cd&\longrightarrow& \cd\\
&((\et1,\ldots,\etd)&,&\vz)&\longmapsto&(\et1 z_1,\ldots, \etd z_d).
\end{array}
$$
This action is effective and Hamiltonian and its moment mapping is
given by
$$
\begin{array}{cccl}
J\,\colon&\cd&\longrightarrow &\rddu\\
&\vz&\longmapsto & \sum_{j=1}^d \zjs
e_j^*+\lambda,\quad\lambda\in\rddu \;\mbox{constant}.
\end{array}
$$
The mapping $J$ is proper and its image is the cone
$\cl=\lambda+\c0$, where $\c0$ denotes the positive orthant in the
space $\rddu$. Now consider the surjective linear mapping
\begin{eqnarray*}
\pi\,\colon &\rd \longrightarrow \rtwo,\\
&e_j \longmapsto X_j
\end{eqnarray*}
and the dimension $2$ quasitorus $D=\rtwo/Q$. Then the linear
mapping $\pi$ induces a quasitorus epimorphism $\Pi\,\colon\,\td
\longrightarrow D$. Define now $N$ to be the kernel of the mapping
$\Pi$ and choose $\lambda=\sum_{j=1}^d \lambda_j e_j^*$. Denote by
$i$ the Lie algebra inclusion $\n=\mbox{Lie}(N)\rightarrow\rd$ and
notice that $\Psi=i^*\circ J$ is a moment mapping for the induced
action of $N$ on $\cd$. Then the quasitorus $\td/N$ acts in a
Hamiltonian fashion on the compact symplectic quasifold
$M=\Psi^{-1}(0)/N$. If we identify the quasitori $D$ and $\td/N$
using the epimorphism $\Pi$, we get a Hamiltonian action of the
quasitorus $D$ whose moment mapping has image equal to
${(\pi^*)}^{-1}(\cl\cap\ker{i^*})=
{(\pi^*)}^{-1}(\cl\cap\mbox{im}\,\pi^*)= {(\pi^*)}^{-1}(\pi^*(\D))$
which is exactly $\D$. This action is effective since the level set
$\Psi^{-1}(0)$ contains points of the form $\vz\in\cd$, $z_j\neq0$,
$j=1,\ldots,d$, where the $T^d$--action is free. Notice finally that
$\dim{M}=2d-2\dim{N}= 2d-2(d-2)=4=2\dim{D}$. If we take $Q$ to be an
ordinary lattice, the space $M$ is either a manifold or an orbifold,
in accordance with the generalization of Delzant's construction to
arbitrary simple rational polytopes by Lerman and Tolman \cite{lt}.

Let us remark that this construction depends on two arbitrary
choices: the choice of the quasilattice $Q$ with respect to which
the polytope is $Q$--rational, and the choice of the
inward--pointing vectors $\xd$ in $Q$.

From now on we fix the quasilattice $Q$ that is generated by the
vectors in the star ${\mathcal S}$ defined by (\ref{star}). It
follows from Subsection~\ref{penrose} that the rhombuses of any
tiling are $Q$--rational with respect to $Q$, in our chosen
coordinate system. The natural choices of inward--pointing vectors
are given by $\pm Y_k,\pm Y_{k+1}$ for $\D^k_R$, and by $\pm Y_k,\pm
Y_{k+2}$ for $\D_r^k$.

Let us begin by performing the generalized Delzant construction for
the thick rhombus $\D_R^2$ and for the thin rhombus $\D_r^4$. We
will show in Theorem~\ref{uguali} that all the other cases can be
reduced to these two.
\subsection{The Thick Rhombus}\label{thick}
\begin{figure}
\begin{center}
\includegraphics{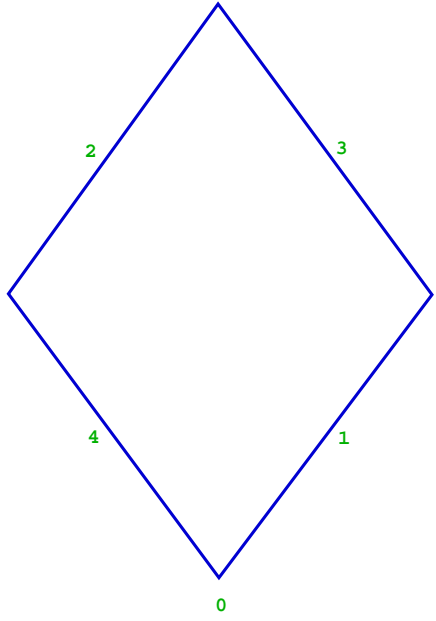}
\end{center}
\caption{The thick rhombus $\D_R^2$} \label{thickrhombus2}
\end{figure}

Let us consider the thick rhombus $\D_R^2$ and let us label its
edges with the numbers $1,2,3,4$, as in Figure~\ref{thickrhombus2}.
The corresponding inward--pointing vectors are given by $X_1=Y_2$,
$X_2=-Y_2$, $X_3=Y_3$, $X_4=-Y_3$, while $\lambda_1=\lambda_4=0$,
and $\lambda_2=\lambda_3= -\frac{1}{2}\sqrt{2+\phi}$. Let us
consider the linear mapping defined by
$$
\begin{array}{cccc}
\pi\colon&\R^4&\rightarrow&\R^2\\
&e_i&\mapsto &X_i.
\end{array}
$$
Its kernel, $\n$, is the $2$-dimensional subspace of $\R^4$ that is
spanned by $e_1+e_2$ and $e_3+e_4$. It is the Lie algebra of
$N=\{\,\exp(X)\in T^4\,|\,X\in\R^4 , \pi(X)\in Q\,\}$. If $\Psi$ is
the moment mapping of the induced $N$--action, then
$$\Psi(z_1,z_2,z_3,z_4)=\left(|z_1|^2+|z_2|^2-\frac{1}{2}\sqrt{2+\phi},
|z_3|^2+|z_4|^2-\frac{1}{2}\sqrt{2+\phi}\right).$$ Let
$R=\left(\frac{1}{2}\sqrt{2+\phi}\right)^{1/2}$ and denote by
$S^3_R$ and $S^2_R$ the spheres of radius $R$, centered at the
origin, of dimension $3$ and $2$ respectively. Then
$\Psi^{-1}(0)=S_R^3\times S_R^3$. A straightforward computation,
using the relations (\ref{rombone}), with $k=2$, gives that
$$N=\left\{\,\exp(X)\in T^4\,|\, X=\left(s,s+\frac{1}{\phi}h,t,t+\frac{1}{\phi}k\right),
s,t \in\R, h,k\in\Z \,\right\},$$ which, for equation~(\ref{phi}),
is equal to
$$
\left\{\,\exp(X)\in T^4\,|\, X=\left(s,s+h\phi,t,t+k\phi\right),
s,t\in\R, h,k\in\Z \,\right\}.
$$
We can think of
\begin{equation}
\label{esse1peresse1} S^1\times S^1=\{\,\exp(X)\in T^4\,|\,
X=(s,s,t,t) , s,t\in\R\,\} \end{equation} as being naturally
embedded in $N$. The quotient group
$$\Gamma=\frac{N}{S^1\times S^1}$$ is a finitely generated group.
In conclusion
$$M_R=\frac{\Psi^{-1}(0)}{N}=\frac{S_R^3\times S_R^3}{N}=
\frac{S^2_R\times S^2_R}{\Gamma}$$ and the quasitorus
$$D^2=\rtwo/Q$$ acts on $M_R$ in a Hamiltonian fashion, with image
of the corresponding moment mapping yielding exactly $\D_R^2$.

It will be useful for the sequel to construct an atlas for the
quasifold $M_R$. It is given by four charts, each of which
corresponds to a vertex of the thick rhombus. Consider for example
the origin, it is given by the intersection of the edges numbered
$1$ and $4$. Let $B_R$ the ball in $\C$ of radius $R$, namely
$$B_R=\{z\in\C\;|\;|z|<R\}.$$
Consider the following mapping, which gives a slice of
$\Psi^{-1}(0)$ transversal to the $N$--orbits
$$\begin{array}{ccc}
B_R\times B_R& \stackrel{\tilde{\tau}_{1,4}}{\longrightarrow}&
\{\vz\in\Psi^{-1}(0)\;|\;z_2\neq0,z_3\neq0\}\\
(z_1,z_4)&\longmapsto&(z_1,\sqrt{R^2-|z_1|^2},\sqrt{R^2-|z_4|^2},z_4)
\end{array}
$$
this induces the homeomorphism
$$
\begin{array}{ccc}
(B_R\times B_R)/\Gamma_{1,4}&\stackrel{\tau_{1,4}}{\longrightarrow}& U_{1,4}\\
\,[\vz]&\longmapsto&[{\tilde{\tau}}_{1,4}(\vz)]
\end{array}
$$
where the open subset $U_{1,4}$ of $M_R$ is the quotient
$$\{\vz\in\Psi^{-1}(0)\;|\;z_2\neq0,z_3\neq0\}/N$$
and the finitely generated group $\Gamma_{1,4}$ is given by
$$\Gamma_{1,4}=N\cap (S^1\times\{1\}\times \{1\}\times S^1)$$
hence
$$\Gamma_{1,4}=\exp\left\{\left(\phi h,0,0,\phi k\right)\;|\; h,k\in\Z\right\}.$$
The triple $(U_{1,4},\tau_{1,4},(B_R\times B_R)/\Gamma_{1,4})$ is a
chart of $M_R$. Analogously we can construct three other charts,
corresponding to the remaining vertices of the thick rhombus, each
of which is characterized by a different pair of variables; the
other three pairs are: $(1,3),(2,3),(2,4)$. These four charts are
compatible, they give therefore an atlas of $M_R$, thus defining on
$M_R$ a quasifold structure.

Now denote by $V_{n}$ the open subset of $S^2_{R}$ given by $S^2_R$
minus the south pole and by $V_{s}$ the open subset of $S^2_R$ given by
$S^2_R$ minus the
north pole, then,  on $\Psi^{-1}(0)$,
consider the action of $S^1\times S^1$ given by
(\ref{esse1peresse1}). We obtain
$$V_n\times V_s=\{\vz\in\Psi^{-1}(0)\;|\;z_2\neq0,z_3\neq0\}/(S^1\times S^1)$$
and
$$U_{1,4}=(V_n\times V_s)/\Gamma.$$
We have the following commutative diagram:
\begin{equation}\label{diagramma}
\xymatrix{
B_R\times B_R\ar[r]^>>>>>>>>{\tilde{\tau}_{1,4}}\ar[d]&\{\vz\in\Psi^{-1}(0)\;|\;z_2\neq0,z_3\neq0\}\ar[d]\\
B_R\times B_R\ar[r]^{\tau_n\times\tau_s}\ar[d]&V_n\times V_s\ar[d]\\
(B_R\times B_R)/\Gamma_{1,4}\ar[r]^{\tau_{1,4}}&U_{1,4} }.
\end{equation}
where the vertical mappings are the natural quotient mappings and
the mappings  $\tau_{n}$ and $\tau_{s}$ from $B_R$ to the open sets
$V_{n}$ and $V_{s}$ respectively are induced by the diagram. Observe
that the mapping
$$\begin{array}{ccc}
\C&\longrightarrow& V_n\\
w&\longmapsto&[\tau_n(R\,w/\sqrt{1+|w|^2})]
\end{array}
$$
is just the stereographic projection from the north pole,
analogously for $\tau_s$. The two charts  $(B_R,\tau_n,V_n)$ and
$(B_R,\tau_s,V_s)$  give a symplectic atlas of $S^2_R$, whose
standard symplectic structure is induced by the standard symplectic
structure on $B_R$. Moreover the symplectic structure of the
quotient $M_R$ is also induced by the standard symplectic structure
on $B_R\times B_R$.

Observe that the quasifold $M_R$ is a global quotient of the product
of spheres by the finitely generated group $\Gamma$; consistently we
have found that its atlas can be obtained by taking the quotient by
$\Gamma$ of the usual atlas of the product of two spheres, given by
the four pairs $V_n\times V_n$, $V_n\times V_s$, $V_s\times V_n$ and
$V_s\times V_s$.
\subsection{The Thin Rhombus}
\begin{figure}
\begin{center}
\includegraphics{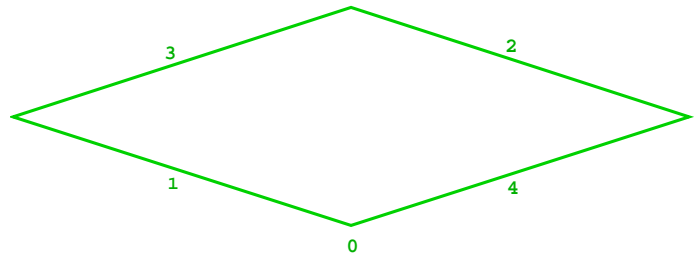}
\end{center}
\caption{The thin rhombus $\D_r^4$} \label{thinrhombus2}
\end{figure}

Let us now consider the thin rhombus $\D_r^4$ and let us label its
edges with the numbers $1,2,3,4$, as in Figure~\ref{thinrhombus2}.
The corresponding inward--pointing vectors are given by $X_1=Y_1$,
$X_2=-Y_1$, $X_3=Y_4$, $X_4=-Y_4$, while $\lambda_1=\lambda_4=0$,
and $\lambda_2=\lambda_3= -\frac{1}{2\phi}\sqrt{2+\phi}$. Let us
consider the linear mapping defined by
$$
\begin{array}{cccc}
\sigma\colon&\R^4&\rightarrow&\R^2\\
&e_i&\mapsto &X_i.
\end{array}
$$
Its kernel, $\l$, is the $2$--dimensional subspace of $\R^4$ that is
spanned by $e_1+e_2$ and $e_3+e_4$. It is the Lie algebra of
$L=\{\,\exp(X)\in T^4\,|\,X\in\R^4 , \sigma(X)\in Q\,\}$. If $\Psi$
is the moment mapping of the induced $L$ action, then
$$\Psi(z_1,z_2,z_3,z_4)=\left(|z_1|^2+|z_2|^2-\frac{1}{2\phi}\sqrt{2+\phi},
|z_3|^2+|z_4|^2-\frac{1}{2\phi}\sqrt{2+\phi}\right).$$ Let
$r=\left(\frac{1}{2\phi}\sqrt{2+\phi}\right)^{1/2}$ and denote by
$S^3_r$ and $S^2_r$ the spheres of radius $r$, centered at the
origin, of dimension $3$ and $2$ respectively. Then
$\Psi^{-1}(0)=S_r^3\times S_r^3$. A straightforward computation,
using the relations (\ref{rombino}), with $k=4$, gives that
$$L=\left\{\,\exp(X)\in
T^4\,|\, X=\left(s,s+h\phi,t,t+k\phi\right), s,t\in\R, h,k\in\Z
\,\right\}=N.$$
In conclusion
$$M_r=\frac{\Psi^{-1}(0)}{N}=\frac{S_r^3\times S_r^3}{N}
=\frac{S_r^2\times S_r^2}{\Gamma}$$ and the quasitorus
$$D^2=\rtwo/Q$$ acts on $M_r$ in a Hamiltonian fashion, with image
of the corresponding moment mapping yielding exactly $\D_r^2$.
Notice that the quasitorus $D^2$ is the same for both the thick and
thin rhombuses.
\section{Symplectic Interpretation of the Tiling} Recall that we
denoted by $M_R$ the symplectic quasifold associated to the thick
rhombus $\D_R^2$ and by $M_r$ the symplectic quasifold associated
with the thin rhombus $\D_r^4$. Consider the five distinguished
thick rhombuses $\D_R^k$ and the five distinguished thin rhombuses
$\D_r^k$, $k=0,\dots,4$. Recall that each of these rhombuses has a
natural choice of inward--pointing vectors, these are $\pm Y_k,\pm
Y_{k+1}$ for $\D^k_R$, and $\pm Y_k,\pm Y_{k+2}$ for $\D_r^k$.
Consider now a Penrose tiling with edges of length $1$. Remark that,
by Proposition~\ref{rotazioni}, in our choice of coordinates, {\em
each} of its rhombuses can be obtained by translation from one of
the $10$ rhombuses $\D_R^k$ and $\D_r^k$. We can then prove the
following
\begin{thm}\label{uguali}The compact symplectic quasifold
corresponding to each thick rhombus of a Penrose tiling with edges
of length $1$ is given by $M_R$. The compact symplectic quasifold
corresponding to each thin rhombus is given by $M_r$.\end{thm}
\proof Observe that, for each $k=0,1,3,4$, there exists a rotation
$P$ of $\R^2$ that leaves the quasilattice $Q$ invariant, that sends
the orthogonal vectors relative to the rhombus $\D_R^k$ to the
orthogonal vectors relative to the rhombus $\D_R^2$, and such that
the dual transformation $P^*$ sends the rhombus $\D_R^2$ to the
rhombus $\D_R^k$. This implies that the reduced space corresponding
to each of the $5$ rhombuses $\D_R^k$, $k=0,\ldots,4$, with the
choice of orthogonal vectors and quasilattice specified above, is
exactly $M_R$. This yields a unique symplectic quasifold, $M_R$, for
all the rhombuses considered. We argue in the same way for each thin
rhombus $\D_r^k$, $k=0,\ldots,4$. Finally, translating the rhombuses
$\D_R^k$ does not produce any change in the corresponding quotient
spaces, therefore, by Proposition~\ref{rotazioni} we are done. \qed

\begin{thm}The quasifolds $M_R$ and $M_r$ are diffeomorphic
but not symplectomorphic.\end{thm} \proof We recall from~\cite{p1}
that a quasifold diffeomorphism, $f$, is a bijective mapping such
that, for each point $p\in M_R$, there is a local model $W$ around
$p$ and a local model $f(W)$ around $f(p)$ such that the mapping
$f$, restricted to $W$, is a diffeomorphism of models as given by
Definition~\ref{diffeomodelli}. A local model is a submodel of a
chart of the atlas that defines the quasifold structure. It is
straightforward to check that, since the manifolds $S^2_R\times
S^2_R$ and $S^2_r\times S^2_r$ are diffeomorphic, the quasifolds
$M_R=(S^2_R\times S^2_R)/\Gamma$ and $M_r=(S^2_r\times
S^2_r)/\Gamma$ are diffeomorphic.

We prove now that $M_R$ and $M_r$ are not symplectomorphic. Denote
by $\omega_R$ and $\omega_r$ the symplectic forms of $M_R$ and $M_r$
respectively. Suppose that there is a symplectomorphism
$h:M_R\longrightarrow M_r$, namely a diffeomorphism $h$ such that
$h^*(\omega_r)=\omega_R$.  The quasifold structures on $M_R$ and
$M_r$ are each defined by four charts, one corresponding to each
vertex of the rhombus, as shown in Subsection~\ref{thick}. We recall
from Remark~\ref{gruppino} that to each point $p\in M_R$ one can
associate finitely generated groups $\Gamma_p$ and $\Lambda_{h(p)}$.
It is easy to check, using Lemma~\ref{green}, that the fact that $h$
is a diffeomorphism implies that these two groups are isomorphic. It
follows from this that $h$ defines a one--to--one correspondence
between the above--given charts of $M_R$ and $M_r$, and that it
sends each of the charts of $M_R$ diffeomorphically onto the
corresponding chart of $M_r$.

Consider now the restriction of $h$ to one such chart of $M_R$, say
$U_{1,4}$, we want to prove that we can construct a diffeomorphism
$\bar{h}$ from $B_R\times B_R$ to $B_r\times B_r$ that lifts the
restriction of $h$ to the given chart. Notice that all submodels of
$U_{1,4}$ are of the type $\wt/\G_{1,4}$, where $\wt$ is an open
subset of $B_R\times B_R$ which, since it is $\G_{1,4}$--invariant,
can either be the product of two open disks, or the product of an
open disk by an open annulus or the product of two open annuli.
Therefore a local model around the point $[\tau_{1,4}(0,0)]$ is
given by $\wt_0=B_{\delta}\times B_{\delta}$, for a suitable
$\delta<R$, modulo the action of $\G_{1,4}$. Since $\wt_0$ is simply
connected here a lift of $h$ is well defined. Moreover, it follows
from the fact that $\phi$ is irrational that when a lift is well
defined on one point of a local model $\wt$ where the action of
$\G_{1,4}$ is free, then it is well defined on all of $\wt$, without
the need of taking a covering of $\wt/\G_{1,4}$. Finally, observe
that if two submodels $\wt_1$ and $\wt_2$ overlap and there is a
lift of $h$ defined on each of them, then by Lemma~\ref{orange},
there is a unique lift defined on $\wt_1\cup \wt_2$. Now the lift
$\bar{h}$ from $B_R\times B_R$ to $B_r\times B_r$ can be constructed
by gluing the local lifts of $h$ that are defined on suitable
submodels of $U_{1,4}$.

Observe now that the four charts of $M_R$ intersect in the
$4$--dimensional dense open subset where the action of the
quasitorus $D^2$ is free; then Lemma~\ref{orange} together with
diagram~(\ref{diagramma}) allow us to lift the diffeomorphism $h$ to
a global diffeomorphism $\tilde{h}$ from $S^2_R\times S^2_R$ to
$S^2_r\times S^2_r$ that is equivariant with respect to the actions
of $\Gamma$ and $\Lambda$ respectively. Moreover, since diagram
(\ref{diagramma}) preserves the symplectic structures, we have that
$\tilde{h}$ is a symplectomorphism between $S^2_R\times S^2_R$ to
$S^2_r\times S^2_r$, which is impossible. \qed

In conclusion there is a {\em unique} quasifold structure that is
naturally associated to the Penrose rhombus tiling, and two distinct
symplectic structures that distinguish the thick and the thin
rhombuses.

\noindent
\small{\sc
Dipartimento di Matematica Applicata "G.
Sansone", Universit\`a di Firenze, Via S. Marta 3, 50139 Firenze,
Italy},\\
{\em E-mail address:} {\tt fiammetta.battaglia@unifi.it}

\medskip

\noindent
\small{\sc Dipartimento di Matematica e Applicazioni per
l'Architettura, Universit\`a di Firenze, Piazza Ghiberti 27, 50122
Firenze, Italy},\\
{\em E-mail address:} {\tt elisa.prato@unifi.it}

\end{document}